\documentclass{conm-p-l}

\usepackage{}

\newtheorem{theorem}{Theorem}[section]

\theoremstyle{definition}

\theoremstyle{remark}

\numberwithin{equation}{section}

\begin{document}

\title{Serre's genus fifty example}


\author{Jaap Top}
\address{Bernoulli Institute, University of Groningen,
Nijenborgh~9, 9747 AG Groningen, the Netherlands.}
\email{j.top@rug.nl}
\thanks{}



\date{}

\begin{abstract}
This note presents explicit equations (up to
birational equivalence over $\mathbb{F}_2$) for
a complete, smooth, absolutely irreducible curve $X$ over
$\mathbb{F}_2$ of genus $50$ satisfying
$\#X(\mathbb{F}_2)=40$. In his 1985 Harvard lecture notes on curves over finite fields, J-P.~Serre
already showed the existence of such a curve:
he used class field theory to describe the
function field $\mathbb{F}_2(X)$ as a certain
abelian extension of the function field
$\mathbb{F}_2(E)$ of some elliptic curve
$E/\mathbb{F}_2$. Although various more
recent texts recall Serre's construction,
explicit equations as well as a description of
intermediate curves $X\to Y\to E$ over $\mathbb{F}_2$ seem to be new. We also describe
explicit equations for a curve over $\mathbb{F}_2$
of genus $8$ with $11$ rational points, and
for a curve over $\mathbb{F}_2$ of genus $22$
with $21$ rational points.
\end{abstract}

\maketitle


\section{Introduction}
Given a smooth and complete, absolutely irreducible
curve $X$ of genus $g$ defined over a finite field
$\mathbb{F}_q$ of cardinality $q$, one can bound
its number of rational points $\#X(\mathbb{F}_q)$
in terms of the two integers $g, q$. As a
consequence, the integer $N_q(g)$, defined as
the maximal $\#X(\mathbb{F}_q)$ that occurs when
$X$ runs over all such curves of genus $g$,
exists. Much of J-P.~Serre's 1985 Harvard lectures
\cite{SerreHarvard} discuss
these integers $N_q(g)$. In particular, the second
part of his lectures considers $g\mapsto N_q(g)$
for fixed prime power $q$, finishing with
a treatment of $N_2(g)$.
One of the results is the estimate
\[N_2(g)\leq 0.6272\cdot g + 9.562 \]
which, in the case $g=50$, shows that any
complete, smooth, absolutely irreducible curve
$X/\mathbb{F}_2$ of genus $50$ satisfies
$\#X(\mathbb{F}_2)\leq 40$.
 
 The last pages of \cite{SerreHarvard} contain
 a proof of the equality $N_2(50)=40$.
 Serre starts from
 an elliptic curve $E/\mathbb{F}_2$ with
 $\#E(\mathbb{F}_2)=5$. He shows using class
 field theory that a covering 
 $\pi\colon X\to E$ over $\mathbb{F}_2$ of degree $8$ exists,
 ramifying only over the Galois orbit of
 some point in $E(\mathbb{F}_{128})\setminus E(\mathbb{F}_2)$, and all $40$ points in
 $\pi^{-1}(E(\mathbb{F}_2))$ are defined over
 $\mathbb{F}_2$, and $X$ has genus $50$.
 
 The aim of this note is to make this explicit.
 Section~\ref{two} briefly describes some properties
 of the elliptic curve $E/\mathbb{F}_2$.
 Section~\ref{three} constructs the desired
 covering $X\to E$ as a tower
 $X=X_{123}\to X_{12}\to X_1\to E$ of covers of degree~$2$.
 The final section discusses a plane model of $X$, and gives remarks on intermediate curves
 $X\to Y\to E$. It ends with an explicit description of curves with the same genus as
 $X_{12}$ and $X_1$ respectively, but with more
 $\mathbb{F}_2$-rational points.
 
\section{The elliptic curve}\label{two}
An elliptic curve over $\mathbb{F}_2$ with
$5=2+1+2$ rational points is supersingular,
hence has $j$-invariant $0$. Up to
isomorphism over $\mathbb{F}_2$ there exist
precisely $3$ elliptic curves over $\mathbb{F}_2$
with $j$-invariant $0$ (see, e.g., \cite[Prop.~3.1]{KST}), with $1,\,3,\,5$
rational points, respectively. So there is
a unique curve $E/\mathbb{F}_2$ with
$\#E(\mathbb{F}_2)=5$. Another way to see this
is by applying \cite[Thm.~4.6]{Sch}. In the
remainder of this paper we give this elliptic curve by
the equation
\[
E\;:\quad y^2+y=x^3+x.
\]
Let $\iota\in\text{End}(E)$ be given by
\[ \iota\colon\quad (x,y)\mapsto (x+1,y+x+1).\]
Note that $\iota^2=[-1]$. Denoting by
$F\in\mbox{End}(E)$ the Frobenius endomorphism
$F(x,y)=(x^2,y^2)$, an easy calculation shows
$F=[-1]+\iota$ and hence 
\[\#E(\mathbb{F}_{2^n})=\deg(1-F^n)=\left\{
\begin{array}{ll}
2^n+1-2^{m+1} & \text{if}\; n=2m\equiv 0\bmod 8;\\
2^n+1+2^{m+1} & \text{if}\; n=2m+1\equiv \pm 1\bmod 8;\\
2^n+1 & \text{if}\;n\equiv \pm 2\bmod 8;\\
2^n+1-2^{m+1} & \text{if}\; n=2m+1\equiv \pm 3\bmod 8;\\
2^n+1+2^{m+1} & \text{if}\; n=2m\equiv 4\bmod 8.
\end{array}
\right.
\]
The case $n=7$ gives $E(\mathbb{F}_{128})=145$,
so $E(\mathbb{F}_{128})\setminus E(\mathbb{F}_2)$
consists of $20$ orbits of points under the action
of $\text{Gal}(\mathbb{F}_{128}/\mathbb{F}_2)$.
One of these orbits is described in the
PhD thesis \cite[\S~9.3]{Auer} of R.~Auer, and
we will use it in the next section. It is
described as follows.

On $E$, consider the function  $y+x^6+x^5+x^2+x$. 
Its only pole is $O=(0:1:0)$, of order $12$. The zeroes of the function
are the points $(a,a^6+a^5+a^2+a)\in E$ where
$a$ satisfies $(a^6+a^5+a^2+a)^2+(a^6+a^5+a^2+a)+a^3+a=0$, i.e.,
\[a^3\cdot(a+1)^2(a^7+a+1)=0.\]
So the zeroes are
$(0,0)$ (with multiplicity $3$), and
$(1,0)$ (with multiplicity $2$), and one
Galois orbit consisting of the $7$ points
$(a,a^6+a^5+a^2+a)\in E(\mathbb{F}_{128})$ where
$a$ runs over the zeroes of $T^7+T+1\in\mathbb{F}_2[T]$.

Write $D$ for the divisor on $E$ consisting
of these $7$ points:
\[
D:=\sum_{a\;\text{such that}\;a^7+a+1=0}(a,a^6+a^5+a^2+a)\in \text{Div}(E).
\]
\section{Constructing the curve via Artin-Schreier extensions}\label{three}
Given are a prime number $p$ and a curve
(absolutely irreducible) $C/\mathbb{F}_p$ with function field $\mathbb{F}_p(C)$. Any $f\in \mathbb{F}_p(C)$
such that the polynomial $T^p-T+f$ is
irreducible in $\overline{\mathbb{F}_p}(C)[T]$
yields a cyclic covering $C_f\to C$ of degree $p$,
namely the curve $C_f/\mathbb{F}_p$ with
function field the (separable, cyclic) extension of $\mathbb{F}_p(C)$
obtained by adjoining the zeros of $T^p-T+f$.
The covering $C_f\to C$ is called an Artin-Schreier covering; its branch points are
the points $P\in C$ that
are poles of $f$ and that have the additional
property that no $g\in \overline{\mathbb{F}_p}(C)$ exists such that $f+g^p-g$ is regular in $P$. 
In particular, all poles of $f$ with order
prime to $p$ will be branch points of $C_f\to C$.
If $P\in C(\mathbb{F}_p)$ is one of the branch points, then the unique point in $C_f$ mapping to it is $\mathbb{F}_p$-rational. If $Q\in C_f(\mathbb{F}_p)$ does not map to a branch point,
then either it lies over a (rational) zero of $f$,
or it lies over a pole $P\in C(\mathbb{F}_p)$ of $f$ and $g$ exists such that $f+g^p-g$ has a zero in $P$.
In particular, over each of the $\mathbb{F}_p$-rational zeros of $f$ one finds
$p$ points in $C_f(\mathbb{F}_p)$.

 The above well known facts will now be used,
 continuing the discussion in Section~\ref{two},
 to construct suitable $C_f\to E$ over $\mathbb{F}_2$ of degree $2$. Since we want
 two rational points over each of the points
 in $E(\mathbb{F}_2)$ we look for $f\in\mathbb{F}_2(E)$ with zeros in all points
 of $E(\mathbb{F}_2)$. To this end, put $R:=\sum P\in\text{Div}(E)$, 
 the divisor consisting of the $5$ rational
 points $P\in E(\mathbb{F}_2)$, and consider
 the Riemann-Roch space $L(D-R)\subset\mathbb{F}_2(E)$ consisting of
 all functions regular away from the points in $D$, with at most simple poles in the points of $D$, and with zeros in the rational points of $E$. The Riemann-Roch theorem implies
 that $L(D-R)$ is a vector space over $\mathbb{F}_2$ of
 dimension $\deg(D-R)=2$. In the present case, 
 $L(D-R)=\{0, f_1, f_2, f_1+f_2\}$ with
 \[ f_1:=\frac{(x^5 + x)\cdot y + (x^2 + x)}{x^7 + x + 1}
 \]
 and
 \[ f_2:=\frac{(x^5 + x^4 + x^3 + x)\cdot y + (x^6 + x^4)}{x^7 + x + 1}.
 \]
 Since the nonzero functions in $L(D-R)$ have a pole of
 order $1$ in the points of $D$, they are not of the form
 $g^2+g$ for any $g\in\overline{\mathbb{F}_2}(E)$.
 Hence adjoining a zero $w_j$ of $T^2+T+f_j$ to $\mathbb{F}_2(E)$,
 one obtains a quadratic extension $\mathbb{F}_2(E)(w_j)$
 which is the function field $\mathbb{F}_2(X_j)$ of
 a smooth and absolutely irreducible curve $X_j$ over $\mathbb{F}_2$. 
 
It turns out that the genus $g(X_j)$ of $X_j$
equals $8$: indeed, consider the differential
$dw_j=df_j$ on $X_j$. Since every pole of $f_j$ in $E$ has order $1$, the differential $df_j$ has 
divisor of poles $2D$. The relation
$w_j(w_j+1)=f_j$ shows that $1/w_j$ is a uniformizer at each of the points of $X_j$ lying
over a point of $D$. The equality $dw_j=w_j^2d(1/w_j)$ then shows that The divisor
of poles of $dw_j$ has degree $14$.
Since $df_j$ (as differential on $E$) has degree
$2g(E)-2=0$, it follows that the divisor of zeros
of $df_j$ (on $E$) has degree $14$. As the degree
two map $X_j\to E$ is unramified away from the poles of $f_j$, the divisor of zeros of $df_j$
on $X_j$ has degree $28$.
As a consequence $2g(X_j)-2=\deg(\text{div}(dw_j))=28-14=14$,
showing that $g(X_j)=8$.
Alternatively, this genus (as well as the
other genera discussed here) can be computed
using the `{\sl F\"{u}hrerdiskriminantenproduktformel}' of
Artin and Hasse; see \cite[VI~\S3 Cor.~2]{SerreCL}
and \cite[Prop.~8.2]{Sch2}. In fact this is how
Serre in \cite{SerreHarvard} finds the relevant
genera.

The polynomial $T^2+T+f_2$ is even irreducible in 
$\overline{\mathbb{F}_2}(X_1)[T]$: it has no zero
in $\overline{\mathbb{F}_2}(E)$, so a zero in
$\overline{\mathbb{F}_2}(X_1)$ would be of
the form $a+bw_1$ for certain 
$a,b\in \overline{\mathbb{F}_2}(E)$ with $b\neq 0$. The form of the zero then shows $b=1$
and $a^2+a=f_1+f_2$. Since $f_1+f_2\in\mathbb{F}_2(E)$ has a pole of exact
order~$1$ in the points of the divisor $D$,
no such $a\in\overline{\mathbb{F}_2}(E)$ exists.
As a consequence the extension 
$\mathbb{F}_2(E)(w_1,w_2)\supset \mathbb{F}_2(E)$
is a Galois extension of degree $4$ (with
group isomorphic to 
$\mathbb{Z}/2\mathbb{Z}\times \mathbb{Z}/2\mathbb{Z}$) 
and it is the function field
of a smooth and absolutely irreducible curve
$X_{12}$ over $\mathbb{F}_2$.

Analogous to the argument above one obtains
$g(X_{12})=22$. Namely, the differential
$dw_2=df_2$ on $X_{12}$ has its only
poles (each of order $2$) at the points over $D$,
so the divisor of poles of $dw_2$ has degree $14$.
The zeros of $df_2$, seen as differential on $E$,
form an effective divisor on $E$ of degree $14$.
Since $X_{12}\to E$ is unramified over the points
in this divisor, one concludes that the divisor
of zeros of $df_2$ on $X_{12}$ has degree
$4\cdot 14=56$. So $2g(X_{12})-2=\deg(\text{div}(dw_j))=56-14=42$,
hence $g(X_{12})=22$.

The next and last step in explicitly constructing
a curve as desired (genus $50$, covering of $E$
with group isomorphic to 
$\mathbb{Z}/2\mathbb{Z}\times \mathbb{Z}/2\mathbb{Z}\times\mathbb{Z}/2\mathbb{Z}$, $40$ points over $\mathbb{F}_2$) is done
by finding a suitable $f_3\in L(2D-R)$.
We demand that $T^2+T+f_3$ is irreducible,
not only as polynomial over
$\overline{\mathbb{F}_2}(E)$ but even in
$\overline{\mathbb{F}_2}(E)(w_1,w_2)[T]$.
To show that $f_3$ exists, substitute $a+bw_1+cw_2+dw_1w_2$ with
$a,b,c,d\in\overline{\mathbb{F}_2}(E)$
into $T^2+T+f_3$.
This leads to the condition
\[
f_3\not\in \wp({\mathbb{F}_2}(E))+L(D-R)
\]
where $\wp\colon {\mathbb{F}_2}(E)\to
{\mathbb{F}_2}(E)$ denotes the
Artin-Schreier map $\xi\mapsto \xi^2+\xi$.
Since we demand $f_3\in L(2D-R)$,
a function $a\in \mathbb{F}_2(E)$
can only satisfy $\wp(a)\in f_3+L(D-R)\subset L(2D-R)$ if $a\in L(D)$. So for $f_3\in L(2D-R)$
to have $T^2+T+f_3$ irreducible over
$\overline{\mathbb{F}_2}(E)(w_1,w_2)$,
a necessary and sufficient condition is that
\[
f_3\not\in \wp(L(D))+L(D-R).
\]

Note that any $a\in L(D)\subset\mathbb{F}_2(E)$
is regular in the points $P\in E(\mathbb{F}_2)$
and $a(P)\in\mathbb{F}_2$. Hence
$\wp(L(D))\subset L(2D-R)$.
The map $\wp\colon L(D)\to L(2D-R)$ is
linear over $\mathbb{F}_2$ with as kernel
the constant functions $\mathbb{F}_2\subset L(D)$.
Moreover $\wp$ maps nonconstant functions in
$L(D)$ to functions in $L(2D-R)$ with a pole
of order $2$ in the points of $D$.
As a consequence $\wp(L(D))\cap L(D-R)=\{0\}$
and
\[
\dim_{\mathbb{F}_2}\left(\wp(L(D))+L(D-R)\right)=(7-1)+2=8.
\]
As $\dim_{\mathbb{F}_2}\,L(2D-R)=\deg(2D-R)=9$,
an $f_3$ as requested exists, and in fact all
$f_3\in L(2D-R)\setminus \left(\wp(L(D))+L(D-R)\right)$
yield the same quadratic extension of
$\mathbb{F}_2(E)(w_1, w_2)$.

One of the possible choices, which we take from now on, is
\[ f_3:=\frac{(x^6+x^5)y+(x^{10}+x^6+x^2+x)}{x^{14}+x^2+1}.
\]
Adjoining a zero $w_3$ of $T^2+T+f_3$ to
$\mathbb{F}_2(E)(w_1,w_2)$, one obtains
the extension $\mathbb{F}_2(E)(w_1,w_2,w_3)\supset
\mathbb{F}_2(E)$ which by construction is Galois
and of degree $8$, 
with Galois group isomorphic to 
$\mathbb{Z}/2\mathbb{Z}\times \mathbb{Z}/2\mathbb{Z}\times\mathbb{Z}/2\mathbb{Z}$, and which is the function field of
a smooth and absolutely irreducible curve
$X_{123}$ over $\mathbb{F}_2$.

The genus of $X_{123}$ can be determined in the
same way as it was done for the other curves
mentioned above:  the function $f_3\in\mathbb{F}_2(E)$ turns out
to have divisor $8\cdot O+3\cdot(1,1)+(0,0)+(0,1)+(1,0)-2D$, hence the
divisor $df_3$ on $E$ has a zero and therefore
also a pole. The only points where it can have
a pole is in the support of $D$, and as $f_3$
has poles of order $2$ here, the pole order of
$df_3$ must necessarily be $2$.
So the degree of the divisor of poles of $df_3$ on $E$ equals $14$. By construction
$X_{123}\to E$ is totally ramified at the points of $D$, so the degree divisor of poles of $df_3$
on $X_{123}$ is also $14$.
The divisor of zeros of $df_3$ on $X_{123}$ has
degree $14\cdot 8=112$, hence on $X_{123}$
one concludes $\deg(\text{div}(df_3))=112-14=98$.
This implies $g(X_{123})=50$.

In fact the differential
$df_3$ on $E$ has divisor $12\cdot O+2\cdot(1,1)-2D$. Since there are $8$ $\mathbb{F}_2$-rational
points on $X_{123}$ over each of the points
in $E(\mathbb{F}_2)$, one obtains
$\#X_{123}(\mathbb{F}_2)=5\cdot 8=40$.
This completes the explicit construction
of a genus $50$ curve $X$ over $\mathbb{F}_2$
with $\#X(\mathbb{F}_2)=40$.

\section{More equations and intermediate curves}\label{four}

For explicit calculations, it is convenient to
present the curve $X$ as a plane curve.
This is done below, starting from the
tower of fields
\[
\mathbb{F}_2(x)\subset \mathbb{F}_2(x,y)
\subset \mathbb{F}_2(x,y,w_1)
\subset  \mathbb{F}_2(x,y,w_1,w_2)
\subset  \mathbb{F}_2(x,y,w_1,w_2,w_3)
\]
in which all successive extensions have degree $2$. The equations are
\[
\left\{
\begin{array}{l}
y^2+y=x^3+x,\\
(x^7+x+1)(w_1^2+w_1)=(x^5+x)\cdot y +x^2+x,\\
(x^7+x+1)(w_2^2+w_2)=(x^5+x^4+x^3+x)\cdot y+x^6+x^4,\\
(x^{14}+x^2+1)(w_3^2+w_3)=(x^6+x^5)\cdot y+
x^{10}+x^6+x^2+x.
\end{array}
\right.
\]
To obtain a plane curve equation, one may start
by observing that the conjugates
$(w_1+a_1)(w_2+a_2)(w_3+a_3)$ (for all
$a_1,a_2,a_3\in\mathbb{F}_2$) of $w_1w_2w_3$
over $\mathbb{F}_2(x,y)$ are pairwise distinct,
hence $\mathbb{F}_2(x,y,w_1,w_2,w_3)=
\mathbb{F}_2(x,y,w_1w_2w_3)$.
A straightforward calculation shows that the
minimal polynomial of $w_1w_2w_3$ over
$\mathbb{F}_2(x,y)$ equals
\[
T^8+T^7+a_6T^6+a_5T^5+a_4T^4+a_5f_1f_2f_3T^3
+a_6f_1^2f_2^2f_3^2T^2+f_1^3f_2^3f_3^3T+f_1^4f_2^4f_3^4
\]
where the $f_j$ are as in Section~\ref{three} and
\[
\begin{array}{c}
a_6:=f_1+f_2+f_3,\\
a_5:=f_1f_2f_3+f_1f_2+f_1f_3+f_2f_3,\\
a_4:=f_1^2f_2^2+f_1^2f_3^2+f_2^2f_3^2+f_1f_2f_3.
\end{array}
\]
This polynomial explicitly defines an extension of
$\mathbb{F}_2(E)=\mathbb{F}_2(x,y)$ as
originally discussed by Serre in \cite{SerreHarvard}.

The polynomial obtained here is not fixed under
the action of $\text{Gal}(\mathbb{F}_2(x,y)/\mathbb{F}_2(x))$
on its coefficients. This is because several
of those coefficients have poles at the points
of the divisor $D$ and are regular in all other
points of $E$; the generator $\gamma$ of
$\text{Gal}(\mathbb{F}_2(x,y)/\mathbb{F}_2(x))$
corresponds to $[-1]\in\text{End}(E)$ and
the support of $[-1]^*D$ is disjoint from that of $D$.
Multiplying the above polynomial by its
conjugate under $\gamma$ results in an
irreducible
degree $16$ element of $\mathbb{F}_2(x)[T]$,
defining $\mathbb{F}_2(x,y,w_1,w_2,w_3)$ as an
extension of $\mathbb{F}_2(x)$ by a single equation.
We do not present it here; it is a
rather long expression in $\mathbb{F}_2[x,\frac{1}{x^7+x+1},T]$, of the form
{\small
\[
T^{16}+\frac{x^{14}+x^{11}+x^{10}+x^6+x^3+x^2+ 1}{(x^7+x+1)^2}T^{14}+
\frac{x^{22} + x^{20} + x^{14} + x^{13} + x^{11} + x^6}{(x^7+x+1)^4}T^{13}+
\ldots
\]
\[
\mbox{ }\hfill{\ldots +
\frac{(x^2+x)^{31}(x^{13} + x^{12} + x^{11} + x^{10} + x^8 + x^5 + x^4 + x + 1)}{(x^7+x+1)^{16}}T +\frac{(x^2+x)^{36}}{(x^7+x+1)^{16}}.}
\]
}
The MAGMA code given below calculates, among other things,
this polynomial {\tt F}.

\vspace{\baselineskip}
\begin{verbatim}
    K0<x>:=FunctionField(GF(2)); P0<Y>:=PolynomialRing(K0);
    KE<y>:=ext<K0 | Y^2+Y+x^3+x>; PE<W>:=PolynomialRing(KE);
    D:=Divisor(Zeros(KE!x^7+x+1)[1]); R:=&+Places(KE,1);
    V,a:=RiemannRochSpace(D-R); RR1:={a(v) : v in V};
    f1:=((x^5+x)*y+x^2+x)/(x^7+x+1); f1 in RR1;
    f2:=((x^5+x^4+x^3+x)*y+x^6+x^4)/(x^7+x+1); f2 in RR1;
    K1<w1>:=ext<KE | W^2+W+f1>; PK1<Z>:=PolynomialRing(K1);
    K12<w2>:=ext<K1 | Z^2+Z+K1!f2>; PK12<T>:=PolynomialRing(K12);
    V2,b:=RiemannRochSpace(2*D-R); RR2:={b(v) : v in V2};
    V3,c:=RiemannRochSpace(D);
    RRred:={a(g)+c(h)+c(h)^2 : g in V , h in V3};
    f3:=((x^6+x^5)*y+x^10+x^6+x^2+x)/(x^14+x^2+1); f3 in RR2;
    f3 in RRred; K123<w3>:=ext<K12 | T^2+T+K12!f3>;
    MinimalPolynomial(w1*w2*w3, KE);
    F:=MinimalPolynomial(w1*w2*w3, K0);
    L<z>:=ext<K0 | F>; Genus(L); #Places(L,1);
\end{verbatim}

Put $S:=\mathbb{F}_2f_1+\mathbb{F}_2f_2+\mathbb{F}_2f_3$. The linear (over $\mathbb{F}_2$) subspaces $R\subset S$ correspond to the fields $F$ with
$\mathbb{F}_2(E)\subset F\subset \mathbb{F}_2(X_{123})$: given $R$, adjoining
the zeros in $\mathbb{F}_2(X_{123})$ of 
$T^2+T+r$ (for all $r\in R$) to $\mathbb{F}_2(E)$
yields the intermediate field $\mathbb{F}_2(X_R)$,
and all intermediate fields are obtained in this way.
The smooth and absolutely irreducible, complete
curves $X_R$ over $\mathbb{F}_2$ satisfy
$g(X_R)=1+7\cdot(-1+\#R)$ and
$\#(X_R)=5\cdot \#R$.
Clearly $X_{\{0\}}=E$ and $X_S=X_{123}$.
Moreover inclusion of subspaces $R\subset R'$ translates into inclusion 
$\mathbb{F}_2(X_R)\subset\mathbb{F}_2(X_{R'})$
of intermediate fields. Is $R\subset R'$, then
$[\mathbb{F}_2(X_{R'}):\mathbb{F}_2(X_{R})]
=\#(R'/R)$.

Although both $E$ and $X_{123}$ reach the
upper bound $N_2(g)$ for the number of
$\mathbb{F}_2$-rational points of curves
of genus $g$ over $\mathbb{F}_2$, 
the `intermediate' curves do not have this
property. If $\dim_{\mathbb{F}_2}R=1$
then $g(X_R)=8$ and
$\#X_R(\mathbb{F}_2)=8<11=N_2(8)$ (for
the latter equality, see
\cite{SerreCR}). And in case
$\dim_{\mathbb{F}_2}R=2$ one has
$g(X_R)=22$ and by \cite[Exc.~8.D]{Sch2} 
$N_2(22)\in\{21,\, 22\}$, whereas $\#X_R(\mathbb{F}_2)=20$.

Quite similar to the construction of the functions
$f_1,f_2,f_3\in\mathbb{F}_2(E)$ and the
resulting explicit curves $X_R$, one finds
equations for curves over $\mathbb{F}_2$
having $N_2(8)=11$ and
$21\leq N_2(22)$ rational points. We briefly sketch how this is done,
and we present the resulting curves.
The function $\frac{x^2+x}{x^3+x+1}$ on
$\mathbb{P}_1$ over $\mathbb{F}_2$ has zeros in
the three points of $\mathbb{P}^1(\mathbb{F}_2)$.
Hence the (hyperelliptic, genus $2$) curve $H$
corresponding to
\[
H\colon y^2+y=\frac{x^2+x}{x^3+x+1}
\]
satisfies $\#H(\mathbb{F}_2)=6$. We will
construct abelian coverings of $H$ with
groups $\mathbb{Z}/2\mathbb{Z}$ and
$\mathbb{Z}/2\mathbb{Z}\times\mathbb{Z}/2\mathbb{Z}$, respectively, in which $5$ points of
$H(\mathbb{F}_2)$ split completely, whereas the
covering is totally ramified over the sixth 
rational point in $H(\mathbb{F}_2)$ and has no other branch points.
As a consequence, the covering has $1+5\cdot d$
rational points, where $d$ is the degree of
the covering map.

As the branch point we take the common zero
$P$ of the two functions $1/x$ and $y(x^3+x+1)/x^3$.
Let $R\in\text{Div}(H)$ be the divisor (of degree
$5$) consisting of the remaining five points
in $H(\mathbb{F}_2)$. The space $L(9P-R)\subset\mathbb{F}_2(H)$ contains the function
\[
a_1:=(x^9+x^7+x^4+x^3)y+x^9+x^8+x^7+x^2
\]
satisfying $\text{ord}_P(a_1)=9$, which implies that
$T^2+T+a_1\in \overline{\mathbb{F}_2}(H)[T]$ is irreducible.
Adjoining a zero of this polynomial to $\mathbb{F}_2(H)$
therefore yields the function field of a smooth, complete
and absolutely irreducible curve $C_1$ over $\mathbb{F}_2$.
It has $g(C_1)=8$ and $\#C_1(\mathbb{F}_2)=11$. Similarly,
$L(11P-R)\subset \mathbb{F}_2(H)$ contains
\[
a_2:=(x^{11}+x^9+x^7+x^6+x^5+x^4+x^3+x^2)y+x^{11}+x^{10}+x^9
+x^7+x^6+x^4
\]
and $\text{ord}_P(a_2)=11$. The polynomial $T^2+T+a_2$ is
not only irreducible over $\overline{\mathbb{F}_2}(H)$ but even
over the quadratic extension 
$\overline{\mathbb{F}_2}(C_1)\supset\overline{\mathbb{F}_2}(H)$.
Hence adjoining a zero of $T^2+T+a_2$ to $\mathbb{F}_2(C_1)$,
an extension is obtained which is Galois
over $\mathbb{F}_2(H)$ with group 
$\mathbb{Z}/2\mathbb{Z}\times\mathbb{Z}/2\mathbb{Z}$, and which is the function field of a smooth, complete and absolutely
irreducible curve $C_2$ over $\mathbb{F}_2$.
It turns out that $g(C_2)=22$ and by construction
$\#C_2(\mathbb{F}_2)=21$. So we obtain the following models, up to birational
equivalence over $\mathbb{F}_2$.
\[
C_1\colon\quad\left\{
\begin{array}{l}
(x^3+x+1)(y^2+y)=x^2+x,\\
z^2+z=a_1
\end{array}
\right.
\]
and
\[
C_2\colon\quad\left\{
\begin{array}{l}
(x^3+x+1)(y^2+y)=x^2+x,\\
z^2+z=a_1,\\
w^2+w=a_2.
\end{array}
\right.
\]

\bibliographystyle{amsplain}
\bibliography{bibliography}

\end{document}